\documentclass[10pt]{article}
\usepackage{graphicx}
\usepackage{amsmath}
\usepackage{amsfonts}
\usepackage{amssymb}

\begin{document}

\title{A NEW TEST FOR THE MULTIVARIATE TWO-SAMPLE PROBLEM BASED ON THE CONCEPT OF
MINIMUM ENERGY}
\author{G. Zech\thanks{Corresponding author. E-mail: zech@physik.uni-siegen.de} \ and
B. Aslan\thanks{E-mail: aslan@physik.uni-siegen.de}\\University of Siegen, Germany}
\maketitle
\begin{abstract}
We introduce a new statistical quantity the \emph{energy }to test whether two
samples originate from the same distributions. The energy is a simple
logarithmic function of the distances of the observations in the variate
space. The distribution of the test statistic is determined by a resampling
method. The power of the energy test in one dimension was studied for a
variety of different test samples and compared to several nonparametric tests.
In two and four dimensions a comparison was performed with the Friedman-Rafsky
and nearest neighbor tests. The two-sample energy test was shown to be
especially powerful in multidimensional applications.
\end{abstract}

\section{INTRODUCTION}

\noindent Let $\mathbf{X}_{1},\mathbf{X}_{2},\ldots,\mathbf{X}_{n}$ and
$\mathbf{Y}_{1},\mathbf{Y}_{2},\ldots,\mathbf{Y}_{m}$ be two samples of
independent random vectors with distributions $F$ and $G$, respectively. The
classical two-sample problem then consists of testing the hypothesis%

\[
H_{0}:F(\mathbf{x})=G(\mathbf{x}),\text{ \ for every }\mathbf{x\in}%
\mathbb{R}^{d},
\]
against the general alternative%

\[
H_{1}:F(\mathbf{x})\neq G(\mathbf{x}),\text{ \ for at least one }\mathbf{x\in
}\mathbb{R}^{d},
\]
where the distribution functions $F$ and $G$ are unknown.

Testing whether two samples are consistent with a single unknown distribution
is a task that occurs in many areas of research. A natural and simple approach
is to compare the first two moments of the sample which measure location and
scale. Many tests of this type can be found in the literature
\cite{Duran,Conover, Büning} but distributions may differ in a more subtle
way. Other tests require binning of data like the power-divergence statistic
test \cite{ReadCressie} and tests of the $\chi^{2}$ type. However, a high
dimensional space is essentially empty, as is expressed in the literature by
the term \emph{curse-of-dimensionality }\cite{Scott}, hence these tests are
rather inefficient unless the sample sizes are large.\ Binning-free tests
based on rank statistics are restricted to univariate distributions, and, when
applied to the marginal distributions, they neglect correlations. The
extension of the Wald-Wolfowitz run test \cite{Wald} and the nearest neighbor
test \cite{Henze} avoid these caveats but it is not obvious that they are
sensitive to all kind of differences in the parent distributions from which
the samples are drawn.

In this paper we propose a new test for the two-sample problem - the
\emph{energy test} - which shows high performance independent of the dimension
of the variate space and which is easy to implement. Our test is related to
Bowman-Foster test \cite{Bowman} but whereas this test is based on probability
density estimation and local comparison, the energy test explores long range correlations.

In Section 2 we define the test statistic $\Phi_{nm}$. Even though the energy
test has been designed for multivariate applications, we apply it in Section 3
to univariate samples because there a unbiased comparison to several well
established univariate tests is easily possible. A selection of examples and
tests investigated by \cite{Büning2} are considered. These are the
Kolmogorov-Smirnov, Cram\`{e}r-von Mises, Wilcox \cite{Wilcox} and Lepage
\cite{Lepage} tests. We have added the $\chi^{2}$ test with equal probability bins.

In Section 4 we study the power of the energy test in two and four dimensions
and compare it to the Friedman-Rafsky \cite{Friedman} and the nearest neighbor tests.

We conclude in Section 5 with a short summary.

\section{THE TWO-SAMPLE \textit{ENERGY }TEST}

The basic idea behind using the quantity \emph{energy} to test the
compatibility of two samples is simple. We consider the sample $A:\mathbf{X}%
_{1},\mathbf{X}_{2},\ldots,\mathbf{X}_{n}$ as a system of positive charges of
charge $1/n$ each, and the second sample $B:\mathbf{Y}_{1},\mathbf{Y}%
_{2},\ldots,\mathbf{Y}_{m}$ as a system of negative charges of charge $-1/m$.
The charges are normalized such that each sample contains a total charge of
one unit. From electrostatics we know that in the limit of where $n,m$ tend to
infinity,$\ $the total potential energy of the combined samples computed for a
potential following a one-over-distance law will be minimum if both charge
samples have the same distribution. The energy test generalizes these
conditions. For the two-sample test we use a logarithmic potential in
$\mathbb{R}^{d}$. In the Appendix we show that also in this case, large values
of energy indicate significant deviations between the parent populations of
the two samples.

\subsection{The test statistic}

The test statistic $\Phi_{nm}$ consists of three terms, which correspond to
the energies of samples $A$ ($\Phi_{A}$), $B$ ($\Phi_{B}$) and the interaction
energy ($\Phi_{AB}$) of the two samples%

\begin{align*}
\Phi_{nm}  & =\Phi_{A}+\Phi_{B}+\Phi_{AB}\\
\Phi_{A}  & =\frac{1}{n^{2}}\sum_{i<j}^{n}R\left(  \left|  \mathbf{x}%
_{i}-\mathbf{x}_{j}\right|  \right)  ,\\
\Phi_{B}  & =\frac{1}{m^{2}}\sum_{i<j}^{m}R\left(  \left|  \mathbf{y}%
_{i}-\mathbf{y}_{j}\right|  \right)  ,\\
\Phi_{AB}  & =-\frac{1}{nm}\sum_{i=1}^{n}\sum_{j=1}^{m}R\left(  \left|
\mathbf{x}_{i}-\mathbf{y}_{j}\right|  \right)  ,
\end{align*}
where $R(r)$ is a continuous, monotonic decreasing function of the Euclidean
distance $r$ between the charges. The choice of $R$ may be adjusted to a
specific statistical problem. For the present analysis, we select $R(r)=-\ln
r$ instead of the electrostatic potential $1/r$. With this choice the test is
scale invariant and offers a good rejection power against many alternatives to
the null hypothesis.

To compute the power of the new two-sample \textit{energy} test we use the
permutation method \cite{Efron} to evaluate the distribution of $\Phi_{nm}$
under $H_{0}$. We merge the $N=m+n$ observations of both samples and draw from
the combined sample a subsample of size $n$ without replacement. The remaining
$m$ observations represent a second sample. The probability distribution under
$H_{0}$ of $\Phi_{nm}$ is evaluated by determining the values of $\Phi_{nm}$
of all $\binom{N}{m}=\frac{N!}{n!m!}$ possible permutations. For large $N$
this procedure can become computationally too laborious. Then the probability
distribution is estimated from a random sample of all possible permutations.

\subsection{Normalization of the distance}

The Euclidean distances between two observations $\mathbf{z}_{i}$ and
$\mathbf{z}_{j}^{\prime}$ in $\mathbb{R}^{d}$ is
\[
\left|  \mathbf{z}_{i}-\mathbf{z}_{j}^{\prime}\right|  =\sqrt{\sum_{k=1}%
^{d}(z_{ik}-z_{jk}^{\prime})^{2}}%
\]
with projections $z_{ik}$ and $z_{jk}^{\prime}$, $k=1,\ldots,d$ of the vectors
$\mathbf{z}_{i}$ and $\mathbf{z}_{j}^{\prime}$.

Since the relative scale of the different variates usually is arbitrary we
propose to normalize the projections by the following transformation
\[
z_{ik}^{^{\ast}}=\frac{z_{ik}-\mu_{k}}{\sigma_{k}}\text{ \ \ \ \ \ \ \ \ }
\begin{array}
[c]{c}%
i=1,\ldots,n\\
k=1,\ldots,d
\end{array}
\]
where $\mu_{k}$, $\sigma_{k}$ are mean value and standard deviation of the
projection $z_{1k},\ldots,z_{nk}$ of the coordinates of the observations of
the pooled sample. In this way we avoid that a single projection dominates the
value of the energy and that other projections contribute only marginally to it.

We did not apply this transformation in the present study, because this might
have biased the comparison with other tests and because the different variates
had similar variances.

\section{POWER COMPARISONS}

The performance of various tests were assessed for finite sample sizes by
Monte Carlo simulations in $d=1$, $2$ and $4$ dimensions. Also the critical
values of all considered tests were calculated by Monte Carlo simulation. We
chose a $5\%$ significance level.

For the null hypothesis we determine the distribution of $\Phi_{nm}$ with the
permutation technique, as mentioned above. We followed \cite{Efron} and
generated $1000$ randomly selected two subsets in each case and determined the
critical values $\phi_{c\text{ }}$of $\phi_{nm}$. For the specific case
$n=m=50$ and samples drawn from a uniform distribution we studied the
statistical fluctuations. Transforming the confidence interval of $\phi_{c}$
into limits for $\alpha$, we obtain the interval $[0.036,0.063]$, see Table
\ref{DistrAlpha}.%

\begin{table}[tbp] \centering
\caption{\label{DistrAlpha}%
Confidence intervals as a function of the number of permutations
for nominal $\alpha=0.05$.}
\begin{tabular}
[c]{cc}\hline
$%
\begin{array}
[c]{c}%
\text{\textbf{\# }of}\\
\text{permutations}%
\end{array}
$ & $CL(95\%)$ for $\alpha$\\\hline
$100$ & $[0.006,0.095]$\\
$300$ & $[0.025,0.075]$\\
$500$ & $[0.031,0.068]$\\
$1000$ & $[0.036,0.063]$\\\hline
\end{tabular}%
\end{table}%

\subsection{One dimensional case}

Even though the energy test has been designed for multivariate applications,
we investigate its power in one dimension because there a comparison with
several well established tests is possible . To avoid a personal bias we drew
the two samples from the probability distributions, which have also been
investigated by \cite{Büning1}:
\begin{align*}
f_{1}(x)  & =\left\{
\begin{array}
[c]{c}%
1\\
0
\end{array}%
\begin{array}
[c]{c}%
-\sqrt{3}\leq x\leq\sqrt{3}\\
\text{otherwise}%
\end{array}
\right. \\
f_{2}(x)  & =\frac{1}{\sqrt{2\pi}}e^{-\frac{x^{2}}{2}}\text{ \ \ \ \ \ \ }\\
f_{3}(x)  & =\frac{1}{2}e^{-\left|  x\right|  }\text{ \ \ \ }\\
f_{4}(x)  & =\frac{1}{\pi}\frac{1}{1+x^{2}}\text{ \ , \ \ \ Cauchy}\\
f_{5}(x)  & =e^{-(x+1)}\text{ \ , \ }\;x\geq-1\text{\ \ \ }\\
f_{6}(x)  & =\chi_{3}^{2}\text{ \ \ \ }
\begin{array}
[c]{c}%
\chi^{2}\text{ with }3\text{ degrees of freedom,}\\
\text{transformed to mean }0\text{, variance }1
\end{array}
\\
f_{7}(x)  & =\frac{1}{2}N(1.5,1)+\frac{1}{2}N(-1.5,1)\\
f_{8}(x)  & =0.8N(0,1)+0.2N(0,4^{2})\\
f_{9}(x)  & =\frac{1}{2}N(1,2^{2})+\frac{1}{2}N(-1,1)
\end{align*}

This set $f_{1}$ to $f_{9}$ of probability distributions covers a variety of
cases of short tailed up to very long tailed probability distributions as well
as skewed ones.

To evaluate the power of the tests we generated $1000$ pairs of samples for
small $n=m=25$, moderate $n=50$, $m=40$ and ``large'' $n=100$, $m=50$, for
seven different scenarios. We have transformed the variates $Y_{i}^{\ast
}=\theta+\tau Y_{j}$, $j=1,\ldots,m$ of the second sample, corresponding to
the alternative distribution, with different location parameters $\theta$ and
scale parameters $\tau$. Powers were simulated in all cases by counting the
number of times a test resulted in a rejection divided by $1000$. All tests
have a nominal significance level of $0.05$.

Table \ref{power1} shows the estimated power for small sample sizes, $n=25$,
$m=25$, of the selected tests. These are the Kolmogorov-Smirnov (KS),
Cram\`{e}r-von Mises (CvM), Wilcox (W), Lepage (L). We have added the
$\chi^{2}$ test with $5$ equal probability bins. Tables \ref{power2} and
\ref{power3} present the results for $n=50$, $m=40$ and $n=100$, $m=50$,
respectively. For the large sample the number of $\chi^{2}$ bins was increased
to $10$.

It is apparent that none of the considered tests performs better than all
other tests for all alternatives. The results indicate that the power of the
energy test in most of the cases is larger than that of the well known
$\chi^{2}$ and KS tests and comparable to that of the CvM test. For long
tailed distributions, e.g. for combinations $(f_{8},f_{4})$, the energy test
is the most powerful test. This is not unexpected since $R(x)=-\ln(x)$ is long
range. Lepage and Wilcox tests are powerful tests for all combinations and
sample sizes considered, however, the Lepage test is based on the first two
moments of the null distribution and therefore specifically adapted to the
type of study presented here.%

\begin{table}[tbp] \centering
\caption{\label{power1} Power of the selected tests for n=m=25,
$\alpha=0.05$, $x\rightarrow\theta+\tau x$ }{}
%
\begin{tabular}
[c]{ccccccccc}\hline
$P_{1}$ & $P_{2}$ & $\theta,\tau$ & KS & CvM & W & L & $\Phi_{25,25}$ &
$\chi^{2}$\\\hline
$f_{1}(x)$ & $f_{7}(x)$ & $%
\begin{array}
[c]{c}%
0.4;1.4\\
0.6,1.6\\
0.6;0.8\\
0.5;0.5
\end{array}
$ & $%
\begin{array}
[c]{c}%
0.12\\
0.37\\
0.40\\
0.70
\end{array}
$ & $%
\begin{array}
[c]{c}%
0.18\\
0.41\\
0.55\\
0.70
\end{array}
$ & $%
\begin{array}
[c]{c}%
0.48\\
0.87\\
0.66\\
0.93
\end{array}
$ & $%
\begin{array}
[c]{c}%
0.38\\
0.69\\
0.50\\
0.86
\end{array}
$ & $%
\begin{array}
[c]{c}%
0.24\\
0.54\\
0.45\\
0.85
\end{array}
$ & $%
\begin{array}
[c]{c}%
0.11\\
0.17\\
0.52\\
0.85
\end{array}
$\\\hline\hline
\multicolumn{1}{l}{$f_{7}(x)$} & \multicolumn{1}{l}{$f_{2}(x)$} &
\multicolumn{1}{l}{$%
\begin{array}
[c]{c}%
0.4;1.4\\
0.6,1.6\\
0.6;0.8\\
0.5;0.5
\end{array}
$} & \multicolumn{1}{l}{$%
\begin{array}
[c]{c}%
0.08\\
0.20\\
0.34\\
0.72
\end{array}
$} & \multicolumn{1}{l}{$%
\begin{array}
[c]{c}%
0.13\\
0.29\\
0.46\\
0.69
\end{array}
$} & \multicolumn{1}{l}{$%
\begin{array}
[c]{c}%
0.22\\
0.46\\
0.57\\
0.93
\end{array}
$} & \multicolumn{1}{l}{$%
\begin{array}
[c]{c}%
0.14\\
0.34\\
0.51\\
0.93
\end{array}
$} & $%
\begin{array}
[c]{c}%
0.13\\
0.31\\
0.44\\
0.89
\end{array}
$ & \multicolumn{1}{l}{$%
\begin{array}
[c]{c}%
0.08\\
0.14\\
0.45\\
0.88
\end{array}
$}\\\hline\hline
\multicolumn{1}{l}{$f_{2}(x)$} & \multicolumn{1}{l}{$f_{3}(x)$} &
\multicolumn{1}{l}{$%
\begin{array}
[c]{c}%
0.4;1.4\\
0.6,1.6\\
0.6;0.8\\
0.5;0.5
\end{array}
$} & \multicolumn{1}{l}{$%
\begin{array}
[c]{c}%
0.17\\
0.33\\
0.64\\
0.74
\end{array}
$} & \multicolumn{1}{l}{$%
\begin{array}
[c]{c}%
0.23\\
0.44\\
0.70\\
0.77
\end{array}
$} & \multicolumn{1}{l}{$%
\begin{array}
[c]{c}%
0.22\\
0.42\\
0.67\\
0.84
\end{array}
$} & \multicolumn{1}{l}{$%
\begin{array}
[c]{c}%
0.19\\
0.37\\
0.66\\
0.91
\end{array}
$} & $%
\begin{array}
[c]{c}%
0.19\\
0.38\\
0.67\\
0.89
\end{array}
$ & \multicolumn{1}{l}{$%
\begin{array}
[c]{c}%
0.14\\
0.24\\
0.60\\
0.84
\end{array}
$}\\\hline\hline
\multicolumn{1}{l}{$f_{2}(x)$} & \multicolumn{1}{l}{$f_{9}(x)$} &
\multicolumn{1}{l}{$%
\begin{array}
[c]{c}%
0.4;1.4\\
0.6,1.6\\
0.6;0.8\\
0.5;0.5
\end{array}
$} & \multicolumn{1}{l}{$%
\begin{array}
[c]{c}%
0.06\\
0.14\\
0.46\\
0.71
\end{array}
$} & \multicolumn{1}{l}{$%
\begin{array}
[c]{c}%
0.09\\
0.22\\
0.59\\
0.72
\end{array}
$} & \multicolumn{1}{l}{$%
\begin{array}
[c]{c}%
0.19\\
0.35\\
0.65\\
0.89
\end{array}
$} & \multicolumn{1}{l}{$%
\begin{array}
[c]{c}%
0.13\\
0.29\\
0.59\\
0.88
\end{array}
$} & $%
\begin{array}
[c]{c}%
0.12\\
0.21\\
0.45\\
0.82
\end{array}
$ & \multicolumn{1}{l}{$%
\begin{array}
[c]{c}%
0.07\\
0.11\\
0.54\\
0.85
\end{array}
$}\\\hline\hline
\multicolumn{1}{l}{$f_{6}(x)$} & \multicolumn{1}{l}{$f_{5}(x)$} &
\multicolumn{1}{l}{$%
\begin{array}
[c]{c}%
0.4;1.4\\
0.6,1.6\\
0.6;0.8\\
0.5;0.5
\end{array}
$} & \multicolumn{1}{l}{$%
\begin{array}
[c]{c}%
0.10\\
0.16\\
0.95\\
0.99
\end{array}
$} & \multicolumn{1}{l}{$%
\begin{array}
[c]{c}%
0.16\\
0.25\\
0.94\\
0.98
\end{array}
$} & \multicolumn{1}{l}{$%
\begin{array}
[c]{c}%
0.15\\
0.24\\
1.00\\
1.00
\end{array}
$} & \multicolumn{1}{l}{$%
\begin{array}
[c]{c}%
0.12\\
0.22\\
0.97\\
1.00
\end{array}
$} & $%
\begin{array}
[c]{c}%
0.12\\
0.16\\
0.98\\
1.00
\end{array}
$ & \multicolumn{1}{l}{$%
\begin{array}
[c]{c}%
0.12\\
0.20\\
0.97\\
1.00
\end{array}
$}\\\hline\hline
\multicolumn{1}{l}{$f_{3}(x)$} & \multicolumn{1}{l}{$f_{8}(x)$} &
\multicolumn{1}{l}{$%
\begin{array}
[c]{c}%
0.4;1.4\\
0.6,1.6\\
0.6;0.8\\
0.5;0.5
\end{array}
$} & \multicolumn{1}{l}{$%
\begin{array}
[c]{c}%
0.26\\
0.55\\
0.85\\
0.90
\end{array}
$} & \multicolumn{1}{l}{$%
\begin{array}
[c]{c}%
0.34\\
0.64\\
0.89\\
0.93
\end{array}
$} & \multicolumn{1}{l}{$%
\begin{array}
[c]{c}%
0.30\\
0.56\\
0.86\\
0.94
\end{array}
$} & \multicolumn{1}{l}{$%
\begin{array}
[c]{c}%
0.25\\
0.53\\
0.84\\
0.97
\end{array}
$} & $%
\begin{array}
[c]{c}%
0.28\\
0.50\\
0.85\\
0.95
\end{array}
$ & \multicolumn{1}{l}{$%
\begin{array}
[c]{c}%
0.21\\
0.43\\
0.79\\
0.91
\end{array}
$}\\\hline\hline
\multicolumn{1}{l}{$f_{8}(x)$} & \multicolumn{1}{l}{$f_{4}(x)$} &
\multicolumn{1}{l}{$%
\begin{array}
[c]{c}%
0.4;1.4\\
0.6,1.6\\
0.6;0.8\\
0.5;0.5
\end{array}
$} & \multicolumn{1}{l}{$%
\begin{array}
[c]{c}%
0.34\\
0.60\\
0.80\\
0.81
\end{array}
$} & \multicolumn{1}{l}{$%
\begin{array}
[c]{c}%
0.42\\
0.67\\
0.85\\
0.84
\end{array}
$} & \multicolumn{1}{l}{$%
\begin{array}
[c]{c}%
0.45\\
0.68\\
0.78\\
0.81
\end{array}
$} & \multicolumn{1}{l}{$%
\begin{array}
[c]{c}%
0.52\\
0.77\\
0.72\\
0.71
\end{array}
$} & $%
\begin{array}
[c]{c}%
0.54\\
0.80\\
0.82\\
0.84
\end{array}
$ & \multicolumn{1}{l}{$%
\begin{array}
[c]{c}%
0.32\\
0.51\\
0.70\\
0.72
\end{array}
$}\\\hline
\end{tabular}%
\end{table}%

\begin{table}[tbp] \centering
\caption{\label{power2} Power of the selected tests for n=50, m=40,
$\alpha=0.05$, $x\rightarrow\theta+\tau x$ }
\begin{tabular}
[c]{ccccccccc}\hline
$P_{1}$ & $P_{2}$ & $\theta;\tau$ & KS & CvM & W & L & $\Phi_{50,40}$ &
$\chi^{2}$\\\hline\hline
$f_{1}(x)$ & $f_{7}(x)$ & $%
\begin{array}
[c]{c}%
0.3;1.3\\
0.4;0.8
\end{array}
$ & $%
\begin{array}
[c]{c}%
0.22\\
0.49
\end{array}
$ & $%
\begin{array}
[c]{c}%
0.18\\
0.47
\end{array}
$ & $%
\begin{array}
[c]{c}%
0.67\\
0.67
\end{array}
$ & $%
\begin{array}
[c]{c}%
0.41\\
0.53
\end{array}
$ & $%
\begin{array}
[c]{c}%
0.25\\
0.46
\end{array}
$ & $%
\begin{array}
[c]{c}%
0.14\\
0.62
\end{array}
$\\\hline\hline
\multicolumn{1}{l}{$f_{7}(x)$} & \multicolumn{1}{l}{$f_{2}(x)$} &
\multicolumn{1}{l}{$%
\begin{array}
[c]{c}%
0.3;1.3\\
0.4;0.8
\end{array}
$} & \multicolumn{1}{l}{$%
\begin{array}
[c]{c}%
0.15\\
0.62
\end{array}
$} & \multicolumn{1}{l}{$%
\begin{array}
[c]{c}%
0.17\\
0.56
\end{array}
$} & \multicolumn{1}{l}{$%
\begin{array}
[c]{c}%
0.34\\
0.66
\end{array}
$} & \multicolumn{1}{l}{$%
\begin{array}
[c]{c}%
0.16\\
0.70
\end{array}
$} & $%
\begin{array}
[c]{c}%
0.20\\
0.58
\end{array}
$ & \multicolumn{1}{l}{$%
\begin{array}
[c]{c}%
0.10\\
0.58
\end{array}
$}\\\hline\hline
\multicolumn{1}{l}{$f_{2}(x)$} & \multicolumn{1}{l}{$f_{3}(x)$} &
\multicolumn{1}{l}{$%
\begin{array}
[c]{c}%
0.3;1.3\\
0.4;0.8
\end{array}
$} & \multicolumn{1}{l}{$%
\begin{array}
[c]{c}%
0.28\\
0.67
\end{array}
$} & \multicolumn{1}{l}{$%
\begin{array}
[c]{c}%
0.29\\
0.66
\end{array}
$} & \multicolumn{1}{l}{$%
\begin{array}
[c]{c}%
0.25\\
0.65
\end{array}
$} & \multicolumn{1}{l}{$%
\begin{array}
[c]{c}%
0.22\\
0.70
\end{array}
$} & $%
\begin{array}
[c]{c}%
0.26\\
0.68
\end{array}
$ & \multicolumn{1}{l}{$%
\begin{array}
[c]{c}%
0.14\\
0.61
\end{array}
$}\\\hline\hline
\multicolumn{1}{l}{$f_{2}(x)$} & \multicolumn{1}{l}{$f_{9}(x)$} &
\multicolumn{1}{l}{$%
\begin{array}
[c]{c}%
0.3;1.3\\
0.4;0.8
\end{array}
$} & \multicolumn{1}{l}{$%
\begin{array}
[c]{c}%
0.07\\
0.51
\end{array}
$} & \multicolumn{1}{l}{$%
\begin{array}
[c]{c}%
0.07\\
0.50
\end{array}
$} & \multicolumn{1}{l}{$%
\begin{array}
[c]{c}%
0.27\\
0.66
\end{array}
$} & \multicolumn{1}{l}{$%
\begin{array}
[c]{c}%
0.12\\
0.58
\end{array}
$} & $%
\begin{array}
[c]{c}%
0.14\\
0.46
\end{array}
$ & \multicolumn{1}{l}{$%
\begin{array}
[c]{c}%
0.07\\
0.61
\end{array}
$}\\\hline\hline
\multicolumn{1}{l}{$f_{6}(x)$} & \multicolumn{1}{l}{$f_{5}(x)$} &
\multicolumn{1}{l}{$%
\begin{array}
[c]{c}%
0.3;1.3\\
0.4;0.8
\end{array}
$} & \multicolumn{1}{l}{$%
\begin{array}
[c]{c}%
0.13\\
1.00
\end{array}
$} & \multicolumn{1}{l}{$%
\begin{array}
[c]{c}%
0.14\\
0.97
\end{array}
$} & \multicolumn{1}{l}{$%
\begin{array}
[c]{c}%
0.18\\
1.00
\end{array}
$} & \multicolumn{1}{l}{$%
\begin{array}
[c]{c}%
0.12\\
0.99
\end{array}
$} & $%
\begin{array}
[c]{c}%
0.11\\
1.00
\end{array}
$ & \multicolumn{1}{l}{$%
\begin{array}
[c]{c}%
0.14\\
1.00
\end{array}
$}\\\hline\hline
\multicolumn{1}{l}{$f_{3}(x)$} & \multicolumn{1}{l}{$f_{8}(x)$} &
\multicolumn{1}{l}{$%
\begin{array}
[c]{c}%
0.3;1.3\\
0.4;0.8
\end{array}
$} & \multicolumn{1}{l}{$%
\begin{array}
[c]{c}%
0.38\\
0.84
\end{array}
$} & \multicolumn{1}{l}{$%
\begin{array}
[c]{c}%
0.39\\
0.84
\end{array}
$} & \multicolumn{1}{l}{$%
\begin{array}
[c]{c}%
0.32\\
0.80
\end{array}
$} & \multicolumn{1}{l}{$%
\begin{array}
[c]{c}%
0.29\\
0.78
\end{array}
$} & $%
\begin{array}
[c]{c}%
0.31\\
0.86
\end{array}
$ & \multicolumn{1}{l}{$%
\begin{array}
[c]{c}%
0.25\\
0.74
\end{array}
$}\\\hline\hline
\multicolumn{1}{l}{$f_{8}(x)$} & \multicolumn{1}{l}{$f_{4}(x)$} &
\multicolumn{1}{l}{$%
\begin{array}
[c]{c}%
0.3;1.3\\
0.4;0.8
\end{array}
$} & \multicolumn{1}{l}{$%
\begin{array}
[c]{c}%
0.50\\
0.76
\end{array}
$} & \multicolumn{1}{l}{$%
\begin{array}
[c]{c}%
0.52\\
0.77
\end{array}
$} & \multicolumn{1}{l}{$%
\begin{array}
[c]{c}%
0.59\\
0.70
\end{array}
$} & \multicolumn{1}{l}{$%
\begin{array}
[c]{c}%
0.66\\
0.63
\end{array}
$} & $%
\begin{array}
[c]{c}%
0.64\\
0.79
\end{array}
$ & \multicolumn{1}{l}{$%
\begin{array}
[c]{c}%
0.39\\
0.59
\end{array}
$}\\\hline
\end{tabular}%
\end{table}%

\begin{table}[tbp] \centering
\caption{\label{power3} Power of the selected tests for n=100, m=50,
$\alpha=0.05$, $x\rightarrow\theta+\tau x$ }
\begin{tabular}
[c]{ccccccccc}\hline
$P_{1}$ & $P_{2}$ & $\theta;\tau$ & KS & CvM & W & L & $\Phi_{100,50}$ &
$\chi^{2}$\\\hline\hline
$f_{1}(x)$ & $f_{7}(x)$ & $%
\begin{array}
[c]{c}%
0.3;1.3\\
0.4;0.8
\end{array}
$ & $%
\begin{array}
[c]{c}%
0.32\\
0.68
\end{array}
$ & $%
\begin{array}
[c]{c}%
0.33\\
0.73
\end{array}
$ & $%
\begin{array}
[c]{c}%
0.97\\
0.85
\end{array}
$ & $%
\begin{array}
[c]{c}%
0.62\\
0.76
\end{array}
$ & $%
\begin{array}
[c]{c}%
0.44\\
0.76
\end{array}
$ & $%
\begin{array}
[c]{c}%
0.28\\
0.66
\end{array}
$\\\hline\hline
\multicolumn{1}{l}{$f_{7}(x)$} & \multicolumn{1}{l}{$f_{2}(x)$} &
\multicolumn{1}{l}{$%
\begin{array}
[c]{c}%
0.3;1.3\\
0.4;0.8
\end{array}
$} & \multicolumn{1}{l}{$%
\begin{array}
[c]{c}%
0.21\\
0.82
\end{array}
$} & \multicolumn{1}{l}{$%
\begin{array}
[c]{c}%
0.27\\
0.79
\end{array}
$} & \multicolumn{1}{l}{$%
\begin{array}
[c]{c}%
0.67\\
0.82
\end{array}
$} & \multicolumn{1}{l}{$%
\begin{array}
[c]{c}%
0.28\\
0.90
\end{array}
$} & $%
\begin{array}
[c]{c}%
0.33\\
0.82
\end{array}
$ & \multicolumn{1}{l}{$%
\begin{array}
[c]{c}%
0.26\\
0.64
\end{array}
$}\\\hline\hline
\multicolumn{1}{l}{$f_{2}(x)$} & \multicolumn{1}{l}{$f_{3}(x)$} &
\multicolumn{1}{l}{$%
\begin{array}
[c]{c}%
0.3;1.3\\
0.4;0.8
\end{array}
$} & \multicolumn{1}{l}{$%
\begin{array}
[c]{c}%
0.31\\
0.85
\end{array}
$} & \multicolumn{1}{l}{$%
\begin{array}
[c]{c}%
0.37\\
0.86
\end{array}
$} & \multicolumn{1}{l}{$%
\begin{array}
[c]{c}%
0.41\\
0.82
\end{array}
$} & \multicolumn{1}{l}{$%
\begin{array}
[c]{c}%
0.29\\
0.90
\end{array}
$} & $%
\begin{array}
[c]{c}%
0.34\\
0.89
\end{array}
$ & \multicolumn{1}{l}{$%
\begin{array}
[c]{c}%
0.17\\
0.72
\end{array}
$}\\\hline\hline
\multicolumn{1}{l}{$f_{2}(x)$} & \multicolumn{1}{l}{$f_{9}(x)$} &
\multicolumn{1}{l}{$%
\begin{array}
[c]{c}%
0.3;1.3\\
0.4;0.8
\end{array}
$} & \multicolumn{1}{l}{$%
\begin{array}
[c]{c}%
0.08\\
0.65
\end{array}
$} & \multicolumn{1}{l}{$%
\begin{array}
[c]{c}%
0.10\\
0.66
\end{array}
$} & \multicolumn{1}{l}{$%
\begin{array}
[c]{c}%
0.51\\
0.79
\end{array}
$} & \multicolumn{1}{l}{$%
\begin{array}
[c]{c}%
0.18\\
0.77
\end{array}
$} & $%
\begin{array}
[c]{c}%
0.21\\
0.67
\end{array}
$ & \multicolumn{1}{l}{$%
\begin{array}
[c]{c}%
0.17\\
0.63
\end{array}
$}\\\hline\hline
\multicolumn{1}{l}{$f_{6}(x)$} & \multicolumn{1}{l}{$f_{5}(x)$} &
\multicolumn{1}{l}{$%
\begin{array}
[c]{c}%
0.3;1.3\\
0.4;0.8
\end{array}
$} & \multicolumn{1}{l}{$%
\begin{array}
[c]{c}%
0.13\\
1.00
\end{array}
$} & \multicolumn{1}{l}{$%
\begin{array}
[c]{c}%
0.19\\
1.00
\end{array}
$} & \multicolumn{1}{l}{$%
\begin{array}
[c]{c}%
0.25\\
1.00
\end{array}
$} & \multicolumn{1}{l}{$%
\begin{array}
[c]{c}%
0.18\\
1.00
\end{array}
$} & $%
\begin{array}
[c]{c}%
0.18\\
1.00
\end{array}
$ & \multicolumn{1}{l}{$%
\begin{array}
[c]{c}%
0.22\\
1.00
\end{array}
$}\\\hline\hline
\multicolumn{1}{l}{$f_{3}(x)$} & \multicolumn{1}{l}{$f_{8}(x)$} &
\multicolumn{1}{l}{$%
\begin{array}
[c]{c}%
0.3;1.3\\
0.4;0.8
\end{array}
$} & \multicolumn{1}{l}{$%
\begin{array}
[c]{c}%
0.46\\
0.95
\end{array}
$} & \multicolumn{1}{l}{$%
\begin{array}
[c]{c}%
0.52\\
0.97
\end{array}
$} & \multicolumn{1}{l}{$%
\begin{array}
[c]{c}%
0.44\\
0.92
\end{array}
$} & \multicolumn{1}{l}{$%
\begin{array}
[c]{c}%
0.44\\
0.94
\end{array}
$} & $%
\begin{array}
[c]{c}%
0.47\\
1.00
\end{array}
$ & \multicolumn{1}{l}{$%
\begin{array}
[c]{c}%
0.25\\
0.86
\end{array}
$}\\\hline\hline
\multicolumn{1}{l}{$f_{8}(x)$} & \multicolumn{1}{l}{$f_{4}(x)$} &
\multicolumn{1}{l}{$%
\begin{array}
[c]{c}%
0.3;1.3\\
0.4;0.8
\end{array}
$} & \multicolumn{1}{l}{$%
\begin{array}
[c]{c}%
0.61\\
0.90
\end{array}
$} & \multicolumn{1}{l}{$%
\begin{array}
[c]{c}%
0.68\\
0.93
\end{array}
$} & \multicolumn{1}{l}{$%
\begin{array}
[c]{c}%
0.81\\
0.87
\end{array}
$} & \multicolumn{1}{l}{$%
\begin{array}
[c]{c}%
0.84\\
0.86
\end{array}
$} & $%
\begin{array}
[c]{c}%
0.86\\
0.92
\end{array}
$ & \multicolumn{1}{l}{$%
\begin{array}
[c]{c}%
0.63\\
0.74
\end{array}
$}\\\hline
\end{tabular}%
\end{table}%

\subsection{Multi-dimensional case}

For the general multivariate two-sample problem, only a few binning-free tests
have been proposed. The Friedman-Rafsky test and the nearest neighbor test
like the energy test are based on the distance between the observations. The
Bowman-Foster goodness-of-fit test uses the probability density estimation
(PDE) to deduce a p.d.f. from a sample. If this technique is applied to both
samples, it obviously can be used as a two sample test. With a Gaussian Kernel
it is almost identical to the energy test with a Gaussian distance function.
We prefer the logarithmic function.

\subsubsection{Friedman-Rafsky test}

The Friedman-Rafsky test can be seen as a generalization of the univariate run
test. The problem in generalizing the run test to more than one dimension is
that there is no unique sorting scheme for the observations. The minimum
spanning tree can be used for this purpose. It is a graph which connects all
observations in such a way that the total Euclidean length of the connections
is minimum. Closed cycles are inhibited.\ The minimum spanning tree of the
pooled sample is formed. The test statistic $R_{nm}$ equals the number of
connections between observations from different samples.

Obviously, in one dimension the test reduces to the run test. Small values of
$R_{nm}\,$lead to a rejection of $H_{0}$. The statistic $R_{nm}$ is
asymptotically distribution-free under the null hypothesis \cite{HenzeRun}.

\subsubsection{The nearest neighbor test}

The nearest neighbor test statistic $N_{nm}$ is the sum of the number of
vectors $\mathbf{Z}_{i}$ of the pooled sample $\left(  \mathbf{Z}_{1}%
,\ldots,\mathbf{Z}_{n+m}\right)  $ where the nearest neighbor of
$\mathbf{Z}_{i}$, denoted by $N(\mathbf{Z}_{i})$, is of the same type as
$\mathbf{Z}_{i}$:
\[
N_{nm}=\sum_{i=1}^{n+m}I\left(  \mathbf{Z}_{i}\text{ and }N(\mathbf{Z}%
_{i})\text{ belong to the same sample}\right)
\]

Here $I$ is the indicator function. $N(\mathbf{Z}_{i})$ can be determined by a
fixed but otherwise arbitrary norm on $\mathbb{R}^{d}$. We selected the
Euclidean norm. In \cite{Henze} it is shown that the limiting distribution of
$N_{nm}$ is normal in the limit $\min(n,m)\rightarrow\infty$ and
$n/(n+m)\rightarrow\tau$ with $0<\tau<1$. Large values of $N_{nm}$ lead to a
rejection of $H_{0}$.

\subsubsection{Comparison of the tests}

\bigskip%
\begin{table}[tbp] \centering
\caption{Two dimensional distributions used to generate the
samples.\label{MultiAlternatives}}
\begin{tabular}
[c]{ccc}\hline
case & $P^{X}$ & $P^{Y}$\\\hline
$1$ & $N\mathbf{(0},\mathbf{I)}$ & $C(\mathbf{0},\mathbf{I)}$\\
$2$ & $N\mathbf{(0},\mathbf{I)}$ & $N_{\log}\mathbf{(}0$,$\mathbf{I)}$\\
$3$ & $N\mathbf{(0},\mathbf{I)}$ & $N\left(  \mathbf{0,}
\begin{array}
[c]{cc}%
1 & 0.6\\
0.6 & 1
\end{array}
\right)  $\\
$4$ & $N\mathbf{(0},\mathbf{I)}$ & $N\left(  \mathbf{0,}
\begin{array}
[c]{cc}%
1 & 0.9\\
0.9 & 1
\end{array}
\right)  $\\
$%
\begin{array}
[c]{c}%
5\\
6
\end{array}
$ & $%
\begin{array}
[c]{c}%
N\mathbf{(0},\mathbf{I)}\\
N\mathbf{(0},\mathbf{I)}%
\end{array}
$ & \ \ \ Student's $\ \
\begin{array}
[c]{c}%
t_{2}\\
t_{4}%
\end{array}
$\\
$%
\begin{array}
[c]{c}%
7\\
8\\
9\\
10\\
11\\
12
\end{array}
$ & $%
\begin{array}
[c]{c}%
U(\mathbf{0},\mathbf{1})\\
U(\mathbf{0},\mathbf{1})\\
U(\mathbf{0},\mathbf{1})\\
U(\mathbf{0},\mathbf{1})\\
U(\mathbf{0},\mathbf{1})\\
U(\mathbf{0},\mathbf{1})
\end{array}
$ & $%
\begin{array}
[c]{c}%
CJ(10)\\
CJ(5)\\
CJ(2)\\
CJ(1)\\
CJ(0.8)\\
CJ(0.6)
\end{array}
$\\
$13$ & $U(\mathbf{0},\mathbf{1})$ & $80\%U(\mathbf{0},\mathbf{1})+20\%N\left(
\mathbf{0.5,}0.05^{2}\mathbf{I}\right)  $\\
$14$ & $U(\mathbf{0},\mathbf{1})$ & $50\%U(\mathbf{0},\mathbf{1})+50\%N\left(
\mathbf{0.5,}0.2^{2}\mathbf{I}\right)  $\\\hline
\end{tabular}%
\end{table}%

In order to investigate how the performance of the tests using $\Phi_{nm}$,
$R_{nm}$ and $N_{nm}$ changes with the dimension, we have considered problems
in dimensions $d=2$ and $4$. In Table \ref{MultiAlternatives} and Table
\ref{MultiAlternatives4} we summarize the alternative probability
distributions $P^{X}$ and $P^{Y}$ from which we drew the two samples. The
first sample was drawn either from $N\mathbf{(0},\mathbf{I)}$ or from
$U(\mathbf{0},\mathbf{1})$ where $N\mathbf{(\mu},\mathbf{V)}$ is a
multivariate normal probability distribution with the indicated mean vector
$\mathbf{\mu}$ and covariance matrix $\mathbf{V}$ and $U(\mathbf{0}%
,\mathbf{1})$ is the multivariate uniform probability distribution in the unit
cube. The parent distributions of the second sample were the Cauchy
distribution $C $, the $N_{\log}$ distribution (explained below), correlated
normal distributions, the Student's distributions $t_{2}$ and $t_{4}$ and
Cook-Johnson $CJ(a)$ distributions \cite{Devroye} with correlation parameter
$a>0$. $CJ(a)$ converges for $a\rightarrow\infty$ to the independent
multivariate uniform distribution and $a\rightarrow0$ corresponds to the
totally correlated case $X_{i1}=X_{i2}=...=X_{id},i=1,...,n$. We generated the
random vectors from $CJ(a)$ via the standard gamma distribution with shape
parameter $a$, following the prescription proposed by \cite{Ahrens}. The
distribution denoted by $N_{\log}$ is obtained by the variable transformation
$x\rightarrow\ln|x|$ applied to each coordinate of a multidimensional normal
distribution and is not to be mixed up with the log-normal distribution. It is
extremely asymmetric. Some of the considered probability densities have also
been used in a power study in \cite{Bahr}.

The various combinations emphasize different types of deviations between the
populations. These include location and scale shifts, differences in skewness
and kurtosis as well as differences in the correlation of the variates.

The test statistics $\Phi_{nm}$, $R_{nm}$ and $N_{nm}$ were evaluated.%

\begin{table}[tbp] \centering
\caption{Four dimensional distributions used to generate the
samples.\label{MultiAlternatives4}}
\begin{tabular}
[c]{ccc}\hline
case & $P^{X}$ & $P^{Y}$\\\hline
$1$ & $N\mathbf{(0},\mathbf{I)}$ & $C(\mathbf{0},\mathbf{I)}$\\
$2$ & $N\mathbf{(0},\mathbf{I)}$ & $N_{\log}\mathbf{(0},\mathbf{I)}$\\
$3$ & $N\mathbf{(0},\mathbf{I)}$ & $80\%N\mathbf{(0},\mathbf{I)+}20\%N\left(
\mathbf{0,}0.2^{2}\mathbf{I}\right)  $\\
$4$ & $N\mathbf{(0},\mathbf{I)}$ & $50\%N\mathbf{(0},\mathbf{I)+}50\%N\left(
\mathbf{0,}
\begin{array}
[c]{cccc}%
1 & 0.4 & 0.5 & 0.7\\
0.4 & 1 & 0.6 & 0.8\\
0.5 & 0.6 & 1 & 0.9\\
0.7 & 0.8 & 0.9 & 1
\end{array}
\right)  $\\
$%
\begin{array}
[c]{c}%
5\\
6
\end{array}
$ & $%
\begin{array}
[c]{c}%
N\mathbf{(0},\mathbf{I)}\\
N\mathbf{(0},\mathbf{I)}%
\end{array}
$ & \ \ \ \ Student's \ \ \ \ \ \ \ \ \ \ $\
\begin{array}
[c]{c}%
t_{2}\\
t_{4}%
\end{array}
$\\
$%
\begin{array}
[c]{c}%
7\\
8\\
9\\
10\\
11\\
12
\end{array}
$ & $%
\begin{array}
[c]{c}%
U(\mathbf{0},\mathbf{1})\\
U(\mathbf{0},\mathbf{1})\\
U(\mathbf{0},\mathbf{1})\\
U(\mathbf{0},\mathbf{1})\\
U(\mathbf{0},\mathbf{1})\\
U(\mathbf{0},\mathbf{1})
\end{array}
$ & $%
\begin{array}
[c]{c}%
CJ(10)\\
CJ(5)\\
CJ(2)\\
CJ(1)\\
CJ(0.8)\\
CJ(0.6)
\end{array}
$\\
$13$ & $U(\mathbf{0},\mathbf{1})$ & $80\%U(\mathbf{0},\mathbf{1}%
)+20\%N(\mathbf{0.5,}0.05^{2}\mathbf{I)}$\\
$14$ & $U(\mathbf{0},\mathbf{1})$ & $50\%U(\mathbf{0},\mathbf{1}%
)+50\%N(\mathbf{0.5,}0.2^{2}\mathbf{I)}$\\\hline
\end{tabular}%
\end{table}%

The power was again computed for $5\%$ significance level \ and samples of
equal size $n=m=30$, $50$, and $100$ (small, moderate and large) in two and
four dimensions. Table \ref{MultiPower2} and Table \ref{MultiPower4}
illustrate the power of the three considered tests calculated from $1000$ replications.%

\begin{table}[tbp] \centering
\caption{Power at significance level $\alpha$ =$ 0.05$, calculated
from $1000$ repetitions, $n=m=30$, $n=m=50$ and $n=m=100$, $d=2$
\label{MultiPower2}}
\begin{tabular}
[c]{ccccccccc}\hline
case &  & $R_{30,30}$ & $N_{30,30}$ & $\Phi_{30,30}$ &  & $R_{50,50}$ &
$N_{50,50}$ & $\Phi_{50,50}$\\\hline
$1$ &  & $0.25$ & $0.23$ & $0.57$ &  & $0.44$ & $0.41$ & $0.86$\\
$2$ &  & $0.33$ & $0.30$ & $0.58$ &  & $0.53$ & $0.50$ & $0.89$\\
$3$ &  & $0.14$ & $0.12$ & $0.13$ &  & $0.17$ & $0.20$ & $0.21$\\
$4$ &  & $0.63$ & $0.57$ & $0.53$ &  & $0.87$ & $0.83$ & $0.87$\\
$5$ &  & $0.14$ & $0.14$ & $0.32$ &  & $0.18$ & $0.20$ & $0.49$\\
$6$ &  & $0.07$ & $0.07$ & $0.11$ &  & $0.08$ & $0.08$ & $0.13$\\
$7$ &  & $0.04$ & $0.07$ & $0.05$ &  & $0.04$ & $0.07$ & $0.05$\\
$8$ &  & $0.05$ & $0.06$ & $0.08$ &  & $0.05$ & $0.08$ & $0.06$\\
$9$ &  & $0.08$ & $0.08$ & $0.10$ &  & $0.08$ & $0.10$ & $0.14$\\
$10$ &  & $0.13$ & $0.12$ & $0.15$ &  & $0.18$ & $0.18$ & $0.23$\\
$11$ &  & $0.15$ & $0.14$ & $0.18$ &  & $0.23$ & $0.22$ & $0.30$\\
$12$ &  & $0.20$ & $0.20$ & $0.25$ &  & $0.33$ & $0.31$ & $0.45$\\
$13$ &  & $0.11$ & $0.10$ & $0.14$ &  & $0.16$ & $0.15$ & $0.33$\\
$14$ &  & $0.09$ & $0.09$ & $0.14$ &  & $0.12$ & $0.11$ & $0.22$\\\hline
\end{tabular}

\bigskip

\bigskip%

\begin{tabular}
[c]{ccccc}\hline
case &  & $R_{100,100}$ & $N_{100,100}$ & $\Phi_{100,100}$\\\hline
$1$ &  & $0.70$ & $0.60$ & $1.00$\\
$2$ &  & $0.82$ & $0.74$ & $1.00$\\
$3$ &  & $0.31$ & $0.28$ & $0.47$\\
$4$ &  & $0.99$ & $0.97$ & $1.00$\\
$5$ &  & $0.34$ & $0.29$ & $0.86$\\
$6$ &  & $0.10$ & $0.11$ & $0.24$\\
$7$ &  & $0.04$ & $0.05$ & $0.10$\\
$8$ &  & $0.05$ & $0.05$ & $0.09$\\
$9$ &  & $0.10$ & $0.11$ & $0.24$\\
$10$ &  & $0.25$ & $0.23$ & $0.52$\\
$11$ &  & $0.32$ & $0.29$ & $0.66$\\
$12$ &  & $0.56$ & $0.48$ & $0.90$\\
$13$ &  & $0.23$ & $0.19$ & $0.78$\\
$14$ &  & $0.16$ & $0.16$ & $0.56$\\\hline
\end{tabular}%
\end{table}%

\begin{table}[tbp] \centering
\caption{Power at significance level $\alpha$ =$ 0.05$, calculated
from $1000$ repetitions, $n=m=30$, $n=m=50$ and $n=m=100$, $d=4$
\label{MultiPower4}}
\begin{tabular}
[c]{ccccccccc}\hline
case &  & $R_{30,30}$ & $N_{30,30}$ & $\Phi_{30,30}$ &  & $R_{50,50}$ &
$N_{50,50}$ & $\Phi_{50,50}$\\\hline
$1$ &  & $0.15$ & $0.19$ & $0.68$ &  & $0.22$ & $0.39$ & $0.93$\\
$2$ &  & $0.46$ & $0.51$ & $0.90$ &  & $0.68$ & $0.78$ & $1.00$\\
$3$ &  & $0.12$ & $0.13$ & $0.23$ &  & $0.14$ & $0.17$ & $0.47$\\
$4$ &  & $0.08$ & $0.17$ & $0.13$ &  & $0.09$ & $0.26$ & $0.18$\\
$5$ &  & $0.16$ & $0.21$ & $0.73$ &  & $0.22$ & $0.35$ & $0.95$\\
$6$ &  & $0.06$ & $0.07$ & $0.17$ &  & $0.08$ & $0.10$ & $0.31$\\
$7$ &  & $0.05$ & $0.06$ & $0.06$ &  & $0.06$ & $0.06$ & $0.07$\\
$8$ &  & $0.06$ & $0.07$ & $0.08$ &  & $0.06$ & $0.09$ & $0.07$\\
$9$ &  & $0.10$ & $0.12$ & $0.18$ &  & $0.13$ & $0.20$ & $0.30$\\
$10$ &  & $0.16$ & $0.26$ & $0.30$ &  & $0.27$ & $0.42$ & $0.62$\\
$11$ &  & $0.23$ & $0.37$ & $0.45$ &  & $0.39$ & $0.58$ & $0.77$\\
$12$ &  & $0.35$ & $0.51$ & $0.65$ &  & $0.58$ & $0.76$ & $0.93$\\
$13$ &  & $0.15$ & $0.16$ & $0.27$ &  & $0.20$ & $0.20$ & $0.62$\\
$14$ &  & $0.11$ & $0.13$ & $0.17$ &  & $0.14$ & $0.18$ & $0.31$\\\hline
\end{tabular}

\bigskip

\bigskip%
\begin{tabular}
[c]{ccccc}\hline
case &  & $R_{100,100}$ & $N_{100,100}$ & $\Phi_{100,100}$\\\hline
$1$ &  & $0.47$ & $0.73$ & $1.00$\\
$2$ &  & $0.93$ & $0.98$ & $1.00$\\
$3$ &  & $0.28$ & $0.25$ & $0.94$\\
$4$ &  & $0.15$ & $0.46$ & $0.49$\\
$5$ &  & $0.45$ & $0.63$ & $1.00$\\
$6$ &  & $0.12$ & $0.16$ & $0.62$\\
$7$ &  & $0.06$ & $0.08$ & $0.10$\\
$8$ &  & $0.05$ & $0.09$ & $0.12$\\
$9$ &  & $0.19$ & $0.29$ & $0.60$\\
$10$ &  & $0.49$ & $0.66$ & $0.97$\\
$11$ &  & $0.69$ & $0.84$ & $0.99$\\
$12$ &  & $0.88$ & $0.96$ & $1.00$\\
$13$ &  & $0.31$ & $0.29$ & $0.99$\\
$14$ &  & $0.22$ & $0.23$ & $0.68$\\\hline
\end{tabular}%
\end{table}%

The \ Friedman-Rafsky and the nearest neighbor tests show very similar
rejection power. For all three sample sizes and dimensions the energy test
performed better than the other two tests in almost all considered
alternatives. This is astonishing because the logarithmic distance function is
long range and the probability distributions in the cases $11$ and $12$ have a
sharp peak in one corner of a $d$ dimensional unit cube and in case $13$ a
sharp peak in the middle of this unit cube. The multivariate student
distribution represents very mild departures from normality, but nevertheless
the rejection rate of the energy test is high.

\section{SUMMARY}

We have introduced the statistic \emph{energy} which is a simple function of
the distances of observations in the sample space. It can be used as a
powerful measure of the compatibility of two samples. It is easy to compute,
efficient and applicable in arbitrary dimensions of the sample space. The
comparison to the Friedman-Rafsky and the nearest neighbor tests demonstrates
its high performance in the multi-dimensional case.

\section{APPENDIX}

We define the energy of the difference of two probability density functions
by
\[
\phi=\frac{1}{2}\int\int\left[  f(\mathbf{x})-f_{0}(\mathbf{x})\right]
\left[  f(\mathbf{x}^{\prime})-f_{0}(\mathbf{x}^{\prime})\right]
R(\mathbf{x},\mathbf{x}^{\prime})d\mathbf{x}d\mathbf{x}^{\prime}.
\]

Here and in what follows, an unspecified integral denotes integration over
$\mathbb{R}^{d}$. Substituting $g(\mathbf{x)=}f(\mathbf{x})-f_{0}(\mathbf{x})$
we obtain
\begin{equation}
\phi=\frac{1}{2}\int\int g(\mathbf{x)}g(\mathbf{x}^{\prime}\mathbf{)}%
R(\mathbf{x},\mathbf{x}^{\prime})d\mathbf{x}d\mathbf{x}^{\prime}%
.\label{energy}%
\end{equation}

We replace in the Eq.(\ref{energy}) the distance function $R(\mathbf{x}%
,\mathbf{x}^{\prime})=R(\left|  \mathbf{x}-\mathbf{x}^{\prime}\right|  )$ by
its Fourier integral
\[
R(\left|  \mathbf{x}-\mathbf{x}^{\prime}\right|  )=\int F(\mathbf{k)}%
e^{i\mathbf{k\cdot(x}-\mathbf{x}^{\prime})}d\mathbf{k}%
\]

and obtain
\begin{align}
\phi & =\frac{1}{2}\int\int\int g(\mathbf{x)}g(\mathbf{x}^{\prime}%
\mathbf{)}F(\mathbf{k)}e^{i\mathbf{k\cdot(x}-\mathbf{x}^{\prime})}%
d\mathbf{x}d\mathbf{x}^{\prime}d\mathbf{k}\nonumber\\
& =\frac{1}{2}\int\left|  G(\mathbf{k)}\right|  ^{2}F(\mathbf{k)}%
d\mathbf{k.}\label{energy1}%
\end{align}

where $G(\mathbf{k)}$ is the Fourier transform of $g(\mathbf{x)}$.

For the function $R\left(  \left|  \mathbf{r}\right|  \right)  =R(r)=\frac
{1}{r^{\kappa}}$ with $d>\kappa$, where $d$ is the dimension of $\mathbf{r}$,
the Fourier transformation $F(\mathbf{k})$ is \cite{Gelfand}:
\[
F(k)=2^{d-\kappa}\pi^{d/2}\frac{\Gamma\left(  \frac{d-\kappa}{2}\right)
}{\Gamma\left(  \frac{\kappa}{2}\right)  }k^{\kappa-d}>0.
\]
with $k=|\mathbf{k}|$.

From Eq.(\ref{energy1}) follows that $\phi$ is positive. The minimum
$\phi_{\min}=0$ is obtained for $g_{\min}(\mathbf{x)\equiv0}$ or
$f(\mathbf{x})\equiv f_{0}(\mathbf{x})$. The result $g_{\min}(\mathbf{x)\equiv
0}$ holds also for the logarithmic distance function $R(r)=-\ln r$ which can
be considered as the $\kappa=0$ limit of the power law distance function:
\[
-\ln r=\underset{n\rightarrow\infty}{\lim}n\left(  \left(  \frac{1}{r}\right)
^{1/n}-1\right)  ,
\]
The additional constant in the distance function does not contribute to the
integral (\ref{energy}).

Expanding (\ref{energy}) we get a sum of three expectation values of
$R(\mathbf{x},\mathbf{x}^{\prime})$
\[
\phi=\frac{1}{2}\int\int\left[  f(\mathbf{x})f(\mathbf{x}^{\prime}%
)-2f_{0}(\mathbf{x})f(\mathbf{x}^{\prime})+f_{0}(\mathbf{x})f_{0}%
(\mathbf{x}^{\prime})\right]  R(\mathbf{x},\mathbf{x}^{\prime})d\mathbf{x}%
d\mathbf{x}^{\prime}%
\]
which can be estimated from two samples drawn from $f$ and $f_{0}$,
respectively, as defined in Section 2:
\begin{multline*}
\phi=\lim_{n,m\rightarrow\infty}\left[  \frac{1}{n(n-1)}\sum_{i<j}^{n}R\left(
\left|  \mathbf{x}_{i}-\mathbf{x}_{j}\right|  \right)  \right.  +\\
-\frac{1}{nm}\sum_{i=1}^{n}\sum_{j=1}^{m}R\left(  \left|  \mathbf{x}%
_{i}-\mathbf{y}_{j}\right|  \right)  +\\
+\left.  \frac{1}{m(m-1)}\sum_{i<j}^{m}R\left(  \left|  \mathbf{y}%
_{i}-\mathbf{y}_{j}\right|  \right)  \right]  .
\end{multline*}

The quantity in the brackets is up to a minor difference in the denominators
equal to our test quantity $\phi_{nm}$. In the limit $n,m\rightarrow\infty$
the statistic $\phi_{nm}$ is minimum if the two samples are from the same distribution.

\end{document}